\documentclass{gtmon_a}
\pdfoutput=1

\proceedingstitle{Proceedings of the School and Conference in Algebraic
Topology (The Vietnam National University, Hanoi, 9-20 August 2004)}
\conferencestart{9 August 2004}
\conferenceend{20 August 2004}
\conferencename{School and Conference in Algebraic Topology}
\conferencelocation{Vietnam National University, Hanoi, Vietnam}

\editor{John Hubbuck}
\givenname{John}
\surname{Hubbuck}

\editor{Nguy\~\ecircumflex{}n H V H\uhorn{}ng}
\givenname{H\uhorn{}ng}
\surname{Nguy\~\ecircumflex{}n}

\editor{Lionel Schwartz}
\givenname{Lionel}
\surname{Schwartz}

\title[Secondary theories for simplicial manifolds]
{Secondary theories for simplicial manifolds\\
and classifying spaces}

\author{Marcello Felisatti}
\givenname{Marcello}
\surname{Felisatti}
\address{Department of Mathematics\\
University of Leicester\\ \newline
Leicester\\ 
LE1 7RH\\
UK}
\email{mf46@mcs.le.ac.uk}
\urladdr{}

\author{Frank Neumann}
\givenname{Frank}
\surname{Neumann}
\email{fn8@mcs.le.ac.uk}
\urladdr{}

\volumenumber{11}
\issuenumber{}
\publicationyear{2007}
\papernumber{3}
\startpage{33}
\endpage{58}

\doi{}
\MR{}
\Zbl{}

\arxivreference{}


\received{24 Mar 2005}
\revised{29 June 2005}
\accepted{11 July 2005}
\published{14 November 2007}
\publishedonline{14 November 2007}
\proposed{}
\seconded{}
\corresponding{}
\editor{}
\version{}

\input xypic
\let\xysavmatrix\xymatrix
\def\xymatrix{\disablesubscriptcorrection\xysavmatrix}
\AtBeginDocument{\let\bar\wbar\let\tilde\wtilde\let\hat\what}

\def\A{{\mathcal A}}

\def\J{{J}}
\def\epsi{\varepsilon}

\def\cc{{\mathcal C}}
\def\d{{\mathrm d}}

\def\NN{{\mathcal N}}

\makeatletter
\def\cnewtheorem#1[#2]#3{\newtheorem{#1}{#3}[section]
\expandafter\let\csname c@#1\endcsname\c@thm}

\theoremstyle{plain}
\newtheorem{thm}{Theorem}[section]
\cnewtheorem{prop}[thm]{Proposition}
\cnewtheorem{cor}[thm]{Corollary}
\cnewtheorem{lem}[thm]{Lemma}

\theoremstyle{definition}
\cnewtheorem{defn}[thm]{Definition}
\cnewtheorem{conj}[thm]{Conjecture}
\cnewtheorem{exmp}[thm]{Example}

\theoremstyle{remark}
\cnewtheorem{rem}[thm]{Remark}
\cnewtheorem{rems}[thm]{Remarks}
\cnewtheorem{note}[thm]{Note}
\cnewtheorem{warn}[thm]{Warning}

\makeatother  %  move after \newtheorem block

\makeop{Mor}
\makeop{op}
\makeop{id}
\makeop{cone}
\begin{document}

\begin{abstract} 
We define secondary theories and characteristic classes for
simplicial smooth manifolds generalizing Karoubi's multiplicative
$K$--theory and multiplicative cohomology groups for smooth manifolds.

As a special case we get versions of the groups of differential
characters of Cheeger and Simons for simplicial smooth manifolds.

Special examples include classifying spaces of Lie groups and Lie
groupoids.
\end{abstract}

\maketitle

\section*{Introduction}

We introduce and analyze secondary theories and characteristic
classes for bundles with connections on simplicial smooth manifolds.

Classical Cheeger--Simons differential characters for simplicial
smooth manifolds with respect to Deligne's `filtration b\^{e}te'
\cite{De} of the associated de Rham complex were first introduced by
Dupont--Hain--Zucker \cite{DHZ} in order to study the relation between the
Cheeger--Chern--Simons invariants of vector bundles with connections on
smooth algebraic varieties and the corresponding characteristic classes
in Deligne--Beilinson cohomology.

In the case of a smooth manifold Dupont, Hain and Zucker showed that
the group of Cheeger--Simons differential characters is isomorphic to
the cohomology group of the cone of the natural map from Deligne's
`filtration b\^{e}te' on the de Rham complex of the manifold to the
complex of smooth singular cochains.

In a series of fundamental papers Karoubi \cite{K1,K2} introduced
multiplicative K--theory and multiplicative cohomology groups, defined
for any filtration of the de Rham complex of a smooth manifold.  By taking
the filtration  to be the `filtration b\^{e}te' it follows that Karoubi's
multiplicative cohomology groups are generalizations of the classical
Cheeger--Simons differential characters in appropriate degrees.

The first author in \cite{F} studied the relationship between
differential characters and multiplicative cohomology further. He
gave a definition of differential characters associated to an
arbitrary filtration of the de Rham complex, which in the case of
the `filtration b\^{e}te' reduces again to the classical case of
Cheeger--Simons. The advantage is that this more general definition
allows for the definition of an explicit map at the levels of
cocycles between Karoubi's multiplicative cohomology groups and
Cheeger--Simons differential characters. It turns out that Karoubi's
multiplicative cohomology groups are the natural gadgets for
systematically constructing and studying secondary characteristic
classes.

Following a similar route in this article we generalize Karoubi's
multiplicative cohomology groups and the groups of Cheeger--Simons
differential characters even further to simplicial smooth manifolds and
arbitrary filtrations of the associated simplicial de Rham complex and
study their relations. This allows for a wider range of applications,
for example to classifying spaces of Lie groups and Lie groupoids.

The outline of the paper is as follows: After introducing the main
background of simplicial de Rham and Chern--Weil theory, mainly
following Dupont \cite{D1,D2} we introduce multiplicative cohomology
groups and groups of differential characters for arbitrary filtrations
of the simplicial de Rham complex. We discuss briefly some examples like
classifying spaces of Lie groups and Lie groupoids. After introducing
the concept of multiplicative bundles  and  multiplicative K--theory
on smooth simplicial manifolds, we construct  characteristic classes of
elements in the multiplicative K--theory with values in multiplicative
cohomology and in the groups of differential characters.

In a sequel to this paper  we will use this approach to construct
and study in a unifying way secondary theories and characteristic
classes for smooth manifolds, foliations, orbifolds, differentiable
stacks etc. basically for everything to which one can associate a
groupoid whose nerve gives rise to a simplicial smooth manifold.
Differential characters for orbifolds were already introduced by
Lupercio and Uribe using closely the approach of Hopkins and Singer
\cite{HS}. Chern--Weil theory for general etale groupoids was
systematically analyzed by Crainic and Moerdijk \cite{CM} using a
very elegant approach to \v{C}ech--de Rham theory, which especially
applies well to leaf spaces of foliated manifolds. Working instead
in the algebraic geometrical context using de Rham theory for
simplicial schemes a similar machinery allows for defining
secondary characteristic classes for Deligne--Mumford stacks, most
prominently for the moduli stack of families of algebraic curves.
Especially multiplicative cohomology with respect to the Hodge or
Hodge--Deligne filtration will be of special interest here. Algebraic
Cheeger--Simons differential characters for algebraic bundles with
connections on smooth algebraic varieties were already studied
systematically by Esnault \cite{E1,E2}.

\subsection*{Acknowledgements}
The first author was supported by EPSRC
research grant GR/S08046/01. The second author likes to thank the
organizers of the Hanoi Conference in Algebraic Topology for
organizing this wonderful conference and for partial financial
support. Final parts of the work were completed while he was
visiting the Tata Institute of Fundamental Research in Mumbai and he
likes to thank the algebraic geometry group and especially 
N\,Nitsure for a very inspiring stay. 
He likes to thank the University of Leicester for granting him a 
study leave.  Both authors also like to thank
J\,R Hunton for many useful discussions.

\section{Elements of simplicial de Rham and Chern--Weil theory}\label{sim}

We recall the ingredients of simplicial de Rham and Chern--Weil theory
as can be found in Dupont \cite{D1,D2} or Dupont--Hain--Zucker \cite{DHZ}.

A simplicial smooth manifold $X_{\bullet}$ is a simplicial object in the
category of ${\mathcal C}^{\infty}$--man\-i\-folds. In other words a simplicial 
smooth manifold is a functor
$$X_{\bullet} \co \Delta^{\op}\rightarrow 
({\mathcal C}^{\infty}-\mathrm{manifolds}).
$$
We can think of $X_{\bullet}$ as a collection $X_{\bullet}=\{X_n\}$ of smooth manifolds $X_n$ for $n\geq 0$ together with smooth face and degeneracy maps
$$\epsi_i \co X_n\rightarrow X_{n-1}, \quad \eta_i \co X_n\rightarrow X_{n+1}$$
for $0\leq i\leq n$ such that the usual simplicial identities hold.
These maps are functorially associated to the inclusion and projection maps
$$\epsi^i \co \Delta^{n-1}\rightarrow \Delta^n, \quad \eta^i \co \Delta^{n+1}\rightarrow \Delta^n.$$
For the differential geometric constructions on $X_{\bullet}$ as introduced below,  the degeneracy maps play no role and everything can be defined for so-called strict simplicial or $\Delta$--manifolds \cite{DHZ}.

The fat realization of a simplicial space $X_{\bullet}$ is the
quotient space
$$\|X_{\bullet}\|= \coprod_{n\geq 0} (\Delta^n \times X_n)/\sim$$
where the equivalence relation is generated by
$$(\epsi^i \times \id)(t,x)\sim (\id \times \epsi_i)(t,x)$$
for any $(t, x) \in \Delta^{n-1} \times X_n$.

There are two versions of the de Rham complex on a simplicial
manifold $X_{\bullet}$ (see Dupont \cite{D1,D2}).

\subsubsection*{The de Rham complex of compatible forms}
A simplicial smooth
complex $k$--form $\omega$ on $X_{\bullet}$ is a sequence
$\{\omega^{(n)}\}$ of smooth complex $k$--forms $\omega^{(n)}\in
\Omega^k_{dR}(\Delta^n\times X_n)$ satisfying the compatibility
condition
$$(\epsi^i\times \id)^* \omega^{(n)}=(\id \times \epsi_i)^*\omega^{(n-1)}$$
in $\Omega^k_{dR}(\Delta^{n-1}\times X_n)$ for all $0\leq i\leq n$
and all $n\geq 1$. Let $\Omega^k_{dR}(X_{\bullet})$ be the set of
all simplicial smooth complex $k$--forms on $X_{\bullet}$. The
exterior differential on $\Omega^k_{dR}(\Delta^n\times X_n)$ induces
an exterior differential $\d$ on $\Omega^k_{dR}(X_{\bullet})$. We
denote by $(\Omega^*_{dR}(X_{\bullet}), \d)$ the de Rham complex of
compatible forms.

Note that $(\Omega^*_{dR}(X_{\bullet}), \d)$ is the total complex of
a double complex $(\Omega^{*,*}_{dR}(X_{\bullet}),\d', \d'')$ with
$$\Omega^k_{dR}(X_{\bullet})=\bigoplus_{r+s=k} \Omega^{r,s}_{dR}(X_{\bullet})$$
and $\d=\d'+ \d''$, where $\Omega^{r,s}_{dR}(X_{\bullet})$ is the
vector space of $(r+s)$--forms, which when restricted to
$\Delta^n\times X_n$ are locally of the form
$$\omega|_{\Delta^n\times X_n}= \sum a_{i_1\ldots i_r j_1\ldots j_s}\d t_{i_1}\wedge
\ldots \wedge \d t_{i_r}\wedge \d x_{j_1}\wedge \d x_{j_s},$$
where $(t_0, \ldots, t_n)$ are barycentric coordinates of $\Delta^n$ and
the $\{x_j\}$ are local coordinates of $X_n$. Furthermore the differentials $\d'$ and $\d''$ are
the exterior differentials on $\Delta^n$ and $X_n$ respectively.

We remark that $\omega=\{\omega^{(n)}\}$ defines a smooth $k$--form on
$$\coprod_{n\geq 0} (\Delta^n \times X_n)$$ and the compatible condition is the necessary and
sufficient condition to define a form on the fat realization
$\|X_{\bullet}\|$ of $X_{\bullet}$
in view of the generating equivalence relation for defining the quotient
space $\|X_{\bullet}\|$.

\subsubsection*{The simplicial de Rham complex}
The de Rham complex
$(\A^*_{dR}(X_{\bullet}), \delta)$ of $X_{\bullet}$ is given as the
total complex of a double complex $(\A^{*,*}_{dR}(X_{\bullet}),
\delta', \delta'')$ with
$$\A^k_{dR}(X_{\bullet})=\bigoplus_{r+s=k} \A^{r,s}_{dR}(X_{\bullet})$$
and $\delta =\delta' + \delta''$, where
$\A^{r,s}_{dR}(X_{\bullet})=\Omega^s_{dR}(X_r)$ is the set of smooth
complex $s$--forms on the smooth manifold $X_r$. Furthermore the
differential
$$\delta''\co \A^{r,s}_{dR}(X_{\bullet})\rightarrow \A^{r, s+1}_{dR}(X_{\bullet})$$
is the exterior differential on $\Omega^*_{dR}(X_r)$ and the
differential
$$\delta'\co\A^{r,s}_{dR}(X_{\bullet})\rightarrow \A^{r+1, s}_{dR}(X_{\bullet})$$
is defined as the alternating sum
$$\delta'=\sum_{i=0}^{r+1}(-1)^i \epsi_i^*.$$

\subsubsection*{The singular cochain complex}
Given a commutative ring $R$ and
a simplicial smooth manifold $X_{\bullet}$ we can also associate a
singular cochain complex $(S^*(X_{\bullet}; R), \partial)$. It is
defined as a double complex $(S^{*,*}(X_{\bullet}; R), \partial',
\partial'')$ with
$$S^k(X_{\bullet}; R)=\bigoplus_{r+s=k} S^{r,s}(X_{\bullet}; R)$$
and $\partial=\partial'+\partial''$, where
$$S^{r,s}(X_{\bullet};R)=S^s(X_r; R)$$
is the set of singular cochains of degree $s$ on the smooth manifold $X_r$.

There is an integration map
$${\mathcal I}\co \A^{r,s}_{dR}(X_{\bullet})\rightarrow S^{r,s}(X_{\bullet}; \C)$$
which gives a morphism of double complexes and Dupont's general version of the de Rham theorem
(see \cite[Proposition 6.1]{D2} for details) shows that this integration map induces natural isomorphisms
$$H(\A^{*}_{dR}(X_{\bullet}, \delta))\cong H(S^*(X_{\bullet},
\C),\partial)\cong H^*(\|X_{\bullet}\|; \C).$$
Stoke's theorem gives  that there is also a morphism of complexes
$${\mathcal J}\co (\Omega^*_{dR}(X_{\bullet}), \d) \rightarrow (\A^*_{dR}(X_{\bullet}), \delta)$$
defined on $\Omega_{dR}^*(\Delta^n\times X_n)$ by integration over
the simplex $\Delta^n$
$$\omega^{(n)}\in \Omega_{dR}^*(\Delta^n\times X_n) \mapsto \int_{\Delta^n} \omega^{(n)}.$$
A result of Dupont \cite[Theorem 2.3 with Corollary 2.8]{D1} gives that this
morphism is in fact a quasi-isomorphism, that is,
$$H(\Omega^*_{dR}(X_{\bullet}), \d)\cong H^*(\A^*_{dR}(X_{\bullet}),
\delta)\cong H^*(\|X_{\bullet}\|, \C).$$

\subsubsection*{The singular cochain complex of compatible cochains}
Let $R$
be a commutative ring. A compatible singular cochain $c$ on
$X_{\bullet}$ is a sequence $\{c^{(n)}\}$ of cochains $c^{(n)}\in
S^k(\Delta^n\times X_n; R)$ satisfying the compatibility condition
$$(\epsi^i\times \id)^* c^{(n)}=(\id \times \epsi_i)^* c^{(n-1)}$$
in $S^k(\Delta^{n-1}\times X_n)$ for all $0\leq i\leq n$ and all
$n\geq 1$. Let $C^k(X_{\bullet}; R)$ be the set of all compatible singular cochains on $X_{\bullet}$ and $(C^*(X_{\bullet}; R), \d)$ be the singular cochain
complex of compatible cochains.

It follows that the natural inclusion of cochain complexes
$$(C^*(X_{\bullet}; R), \d)\rightarrow (S^*(X_{\bullet}; R), \partial)$$
is a quasi-isomorphism (see Dupont--Hain--Zucker \cite{DHZ}).

Integrating forms preserves the compatibility conditions and we
therefore get an induced map of complexes \cite{DHZ}
$${\mathcal I'}\co \Omega^*_{dR}(X_{\bullet})\rightarrow C^*(X_{\bullet}; \C)$$
fitting into a commutative diagram
\[
\xymatrix{ \Omega^*_{dR}(X_{\bullet})\ar[r]^{\mathcal
I'}\ar[d]^{\mathcal J}&
C^*(X_{\bullet}; \C)\ar[d]\\
A^*_{dR}(X_{\bullet})\ar[r]^{\mathcal I}& S^*(X_{\bullet}; \C) }
\]
and which is again a quasi-isomorphism, that is, we have
$$H^*(\Omega^*_{dR}(X_{\bullet}), \d)\cong H^*(C^*(X_{\bullet}; \C), \partial).$$
We will use these compatible de Rham and cochain complexes for the definition of
multiplicative cohomology and differential characters of $X_{\bullet}$ in \fullref{MH}.

We recall finally the basic aspects of Chern--Weil theory in the
simplicial context as developed by Dupont \cite{D2}, and by Dupont,
Hain and Zucker \cite{DHZ}.

\subsubsection*{Principal bundles}
Let $G$ be a Lie group. A principal $G$--bundle over a
simplicial smooth manifold $X_{\bullet}$ is given by a simplicial smooth manifold $P_{\bullet}$ and a morphism $\pi_{\bullet}\co P_{\bullet}\rightarrow X_{\bullet}$ of simplicial smooth manifolds, such that
\begin{itemize}
\item[(i)] for each $n$ the map $\pi_{p}\co P_{n}\rightarrow X_{n}$ is a principal $G$-bundle
over $X_{n}$
\item[(ii)] for each morphism $f\co [m]\rightarrow [n]$ of the simplex category $\Delta$ the induced map
$f^*\co P_{n}\rightarrow P_{m}$ is a morphism of $G$--bundles, that is, we have a commutative
diagram
\[
\xy
\xymatrix{
P_n\ar[r]^{f^*}
\ar[d]&
P_m\ar[d]\\
X_n\ar[r]^{f^*}&X_m}
\endxy
\]
\end{itemize}

It follows, that if $\pi_{\bullet}\co P_{\bullet}\rightarrow X_{\bullet}$ is a principal
$G$--bundle over $X_{\bullet}$, then $|\pi_{\bullet}|\co |P_{\bullet}|\rightarrow |X_{\bullet}|$ is a principal $G$--bundle with $G$--action induced by
$$\Delta^n\times P_{n}\times G\rightarrow \Delta^n\times P_{n}, \, \, (t,x,g)\mapsto (t, xg).$$

\subsubsection*{Connections and curvature on principal bundles}
A connection $\theta$ on a principal $G$--bundle
$\pi_{\bullet}\co P_{\bullet}\rightarrow X_{\bullet}$ over a simplicial manifold
 $X_{\bullet}$ is a $G$--invariant $1$--form (in the de Rham complex of compatible forms)
$$\theta \in \Omega^1_{dR}(P_{\bullet}; \mathfrak{g})$$
taking values in the Lie algebra $\mathfrak{g}$ of $G$,
 on which $G$ acts via the adjoint representation,
 such that for each $n$ the restriction
$$ \theta^{(n)}=\theta|_{\Delta^n\times P_n},$$
is a connection on the bundle
$\pi_n\co \Delta^n\times P_n \rightarrow \Delta^n\times X_n$.
So $\theta=\{\theta^{(n)}\}$ can as well be interpreted as a sequence of
$\mathfrak g$--valued compatible 1--forms.

The curvature $\Omega$ of the connection form $\theta$ is the differential form
$$\Omega=d\theta + \tfrac{1}{2}[\theta, \theta]\in \Omega^2_{dR}(X_{\bullet}; \mathfrak{g}).$$
We have the following general theorem concerning the Chern--Weil map of
a simplicial smooth manifold.

\begin{thm}[Dupont {\cite[Proposition 3.7]{D1}}]\label{CW}
Let $\Phi$ be an invariant polynomial. The differential form
$\Phi(\theta)\in \Omega^*_{dR}(P_{\bullet})$ is a closed form and
descends to a closed form in $\Omega^*_{dR}(X_{\bullet})$ and its
cohomology class represents the image of the class $\Phi\in H^*(BG;
\C)$ under the Chern--Weil map
$$H^*(BG; \C) \rightarrow H^*(\|X_{\bullet}\|; \C)$$ associated to the principal
bundle $\pi_{\bullet}\co P_{\bullet}\rightarrow X_{\bullet}$.
\end{thm}

With an abuse of notation, in the sequel we will denote also by
$\Phi(\theta)$ the form in $\Omega^*_{dR}(X_{\bullet})$.

In order to classify differential geometric invariants on simplicial smooth
manifolds it is useful to extend the constructions outlined above to the category of bisimplicial smooth manifolds. This is straightforward and we will only
briefly describe the constructions (compare also Dupont--Hain--Zucker
\cite{DHZ} and Dupont--Just \cite{DJ}).

A bisimplicial smooth manifold $X_{\bullet \bullet}$ is a 
simplicial object in
the category of simplicial smooth manifolds. In other words 
a bisimplicial smooth manifold is a functor
$$X_{\bullet \bullet}\co \Delta^{\op}\times \Delta^{\op}\rightarrow
({\mathcal C}^{\infty}-\text{manifolds}).$$
We can think of $X_{\bullet \bullet}$ as a collection
$X_{\bullet \bullet}=\{X_{m,n}\}$ of smooth manifolds 
$X_{m,n}$ for $m, n\geq 0$ together with smooth horizontal and 
vertical face and degeneracy maps
$$\epsi_i'\co X_{m,n}\rightarrow X_{m-1, n},\quad 
\epsi_j''\co X_{m,n}\rightarrow X_{m, n-1}$$
$$\eta_i'\co X_{m,n}\rightarrow X_{m+1, n}, \quad 
\eta_j''\co X_{m,n}\rightarrow X_{m,n+1}$$
for $0\leq i\leq m$ and $0\leq j\leq n$, where the horizontal and
vertical maps commute and the usual simplicial identities hold
horizontally and vertically.

The fat realization of a bisimplicial space $X_{\bullet \bullet}$ is
the quotient space
$$\|X_{\bullet \bullet}\|=\coprod_{m,n\geq 0}(\Delta^m\times \Delta^n \times X_{m,n})/\sim$$
where the equivalence relation is generated by
$$(\epsi^i\times \id\times \id)(t,s,x)\sim (\id\times \id \times \epsi'_i)(t,s,x)$$
for any $(t,s,x)\in \Delta^{m-1}\times \Delta^n\times X_{m,n}$ and
$$(\id\times \epsi^j \times \id)(t,s,x)\sim (\id\times \id \times \epsi_j'')(t,s,x)$$
for any $(t,s,x)\in \Delta^m\times \Delta^{n-1}\times X_{m,n}$.

In a similar manner as for simplicial smooth manifolds, we can associate two de
Rham complexes and a singular cochain complex for bisimplicial smooth manifolds.

\subsubsection*{The de Rham complex of compatible forms}
A bisimplicial smooth
$k$--form $\omega$ on $X_{\bullet \bullet}$ is a sequence
$\{\omega^{(m,n)}\}$ of smooth complex $k$--forms
$$\omega^{(m,n)}\in \Omega^k_{dR}(\Delta^m\times \Delta^n \times X_{m,n})$$
satisfying the compatibility conditions
$$(\epsi^i\times \id\times \id)^*\omega^{(m,n)}=(\id\times \id\times \epsi_i')^*\omega^{(m-1,n)}$$
in $\Omega^k_{dR}(\Delta^{m-1}\times \Delta^n \times X_{m,n})$ for
all $0\leq i\leq m$, $m\geq 1$ and $n\geq 0$ as well as the
compatibility conditions
$$(\id\times\epsi^j\times \id)^*\omega^{(m,n)}=(\id\times \id\times \epsi_j'')^*\omega^{(m-1,n)}$$
in $\Omega^k_{dR}(\Delta^m\times \Delta^{n-1} \times X_{m,n})$ for
all $0\leq j\leq n$, $n\geq 1$ and $m\geq 0$.

We denote the set of bisimplicial smooth $k$--forms by
$\Omega^k_{dR}(X_{\bullet \bullet})$. The exterior differential on
$\Omega^*_{dR}(\Delta^m\times \Delta^n \times X_{m,n})$ induces an
exterior differential $\d$ on $\Omega^*_{dR}(X_{\bullet \bullet})$
and we get a complex $(\Omega^*_{dR}(X_{\bullet \bullet}), \d)$, the
de Rham complex of compatible forms on $X_{\bullet \bullet}$.

We note also that we can view the complex $(\Omega^*_{dR}(X_{\bullet
\bullet}), \d)$ as a triple complex
$$(\Omega^{*,*,*}_{dR}(X_{\bullet \bullet}), \d'_{\Delta}, \d''_{\Delta},
\d_X)$$
with
$$\Omega^k_{dR}(X_{\bullet \bullet})=\bigoplus_{r+s+t=k}\Omega^{r,s,t}_{dR}(X_{\bullet \bullet})$$
and $\d=\d'_{\Delta}+ \d''_{\Delta} + \d_X$ where
$\Omega^{r,s,t}_{dR}(X_{\bullet \bullet})$ is the complex vector
space of $(r+s+t)$--forms, which when restricted to $\Delta^m\times
\Delta^n\times X_{m,n}$ are locally of the form
\begin{equation*}
\begin{split}a|_{\Delta^m\times \Delta^n\times X_{m,n}} =
\sum a_{i_1\ldots i_r j_1\ldots j_s k_1\ldots k_t}\d t_{i_1}\wedge\ldots\wedge \d t_{i_r}\wedge & \d s_{j_1}\wedge\ldots\\ \ldots\wedge \d s_{j_s}\wedge \d x_{k_1}\wedge\ldots\wedge \d x_{k_t}
\end{split}
\end{equation*}
with $(t_0,\ldots ,t_m)$ and $(s_0, \dots ,s_n)$ the barycentric coordinates of $\Delta^m$ and $\Delta^n$ respectively and the $\{x_k\}$ are local coordinates
of $X_{m,n}$.

\subsubsection*{The simplicial de Rham complex}
Again we also have the
simplicial de Rham complex $(\A^*(X_{\bullet \bullet}), \delta)$ of
$X_{\bullet \bullet}$ given as the total complex of the triple
complex $$(\A^{*,*,*}_{dR}(X_{\bullet \bullet}), \delta', \delta'',
\delta''')$$ with
$$\A^k_{dR}(X_{\bullet \bullet})=\bigoplus_{r+s+t=k}\A^{r,s,t}_{dR}(X_{\bullet \bullet})$$
with
$$\A^{r,s,t}_{dR}(X_{\bullet \bullet})=\Omega^t_{dR}(X_{r,s})$$
and $\delta=\delta'+\delta''+\delta'''$.

\subsubsection*{The singular cochain complex}
For a commutative ring $R$, we
similarly define the singular cochain complex $(S^*(X_{\bullet
\bullet};R), \partial)$ of $X_{\bullet \bullet}$ given as the total
complex of the triple complex $(S^{*,*,*}(X_{\bullet \bullet}; R),
\partial', \partial'', \partial''')$ with
$$S^k (X_{\bullet \bullet}; R)=\bigoplus_{r+s+t=k}S^{r,s,t}(X_{\bullet \bullet}; R)$$
with
$$S^{r,s,t}(X_{\bullet \bullet}; R)=S^t(X_{r,s}; R)$$
and $\partial=\partial'+\partial''+\partial'''$.

Using iteratively the arguments as in the case for simplicial smooth manifolds, we can finally also derive a de Rham theorem relating the cohomology of all the
complexes defined with the cohomology of the realization of $X_{\bullet
\bullet}$, that is, we have natural isomorphisms
$$H(\Omega^*_{dR}(X_{\bullet \bullet}), \d)\cong H^*(\A^*_{dR}(X_{\bullet
\bullet}), \delta)\cong H^*(\|X_{\bullet \bullet}\|, \C).$$
We remark that we can also define again the singular cochain complex of compatible forms $C^*(X_{\bullet \bullet}; R)$ in a same way as for
$X_{\bullet}$ using two compatibility conditions instead. Again we have
quasi-isomorphisms as in the simplicial case between the various complexes.

Finally we can extend the elements of simplicial Chern--Weil theory to bisimplicial
smooth manifolds, especially we remark that we can define principal $G$--bundles
$$\pi_{\bullet \bullet}\co P_{\bullet \bullet}\rightarrow X_{\bullet \bullet}$$
for the action of a Lie group $G$ and a connection $\theta$ on $\pi_{\bullet \bullet}$ which is again a 1--form
$$\theta \in \Omega^1_{dR}(P_{\bullet \bullet}; \mathfrak{g}).$$
The curvature $\Omega$ of the connection form $\nabla$ is again the
differential form
$$\Omega=d\theta + \tfrac{1}{2}[\theta, \theta]\in \Omega^2_{dR}(X_{\bullet \bullet}; \mathfrak{g}).$$
Again a version of Dupont's theorem (\fullref{CW}) holds in the context
of bisimplicial manifolds. When defining characteristic classes we will
need that given any connection on a principal bundle, we can construct
a connection on (a model of) the universal bundle that pulls back to
the given one. For the convenience of the reader, we recall the  theorem
stating this fact and outline its proof, which for $GL_n(\C)$--principal
bundles is \cite[Proposition 6.15]{DHZ}.

\begin{thm}\label{uni}
Let $G$ be a Lie group, $X_{\bullet}$ a simplicial smooth manifold and $\pi_{\bullet}\co P_{\bullet}\rightarrow X_{\bullet}$ a principal $G$--bundle with connection
$$\theta \in \Omega^1_{dR}(P_{\bullet}; \mathfrak{g}).$$
Then there exists a bisimplicial smooth manifold $B_{\bullet
\bullet}$ of the homotopy type of the classifying space $BG$ and a
$G$--principle bundle $U_{\bullet \bullet}\rightarrow B_{\bullet
\bullet}$ with a connection $\theta_{U_{\bullet \bullet}}\in
\Omega^1_{dR}(U_{\bullet \bullet}; \mathfrak{g})$ and a morphism $(\Psi,
\psi)$ of $G$--bundles
\[
\xy
\xymatrix{
P_{\bullet}\ar[r]^{\Psi}
\ar[d]&U_{\bullet \bullet}\ar[d]\\
X_{\bullet}\ar[r]^{\psi}&B_{\bullet \bullet}}
\endxy
\]
such that $\Psi^*(\theta_{U_{\bullet \bullet}})=\theta$.
\end{thm}

\begin{proof}
We define the bisimplicial manifold $U_{\bullet \bullet}$ as follows:
$$U_{\bullet m}=(P_\bullet)^{m+1}$$
with face maps
\begin{align*}
d_i\co U_{\bullet m}&\longrightarrow U_{\bullet m-1} \\
(u_0, \ldots u_m)&\longmapsto (u_0,\ldots, u_{i-1}, u_{i+1}, \ldots
u_{m})\quad \text{for } 0\leq i\leq m
\end{align*}
and degeneracy maps
\begin{align*}
s_i\co  U_{\bullet m}&\longrightarrow U_{\bullet m+1} \\
(u_0, \ldots u_m)&\longmapsto (u_0,\ldots, u_{i-1},u_i, u_i, u_{i+1}
\ldots u_{m})\quad \text{for } 0\leq i\leq m.
\end{align*}
The fat realization $\|U_{\bullet \bullet}\|$ of this simplicial manifold
is contractible, that is, homotopy equivalent to a point (see Segal
\cite{S}). Now the free $G$-action on $P_{\bullet}$ induces a free
$G$--action on $U_{\bullet \bullet}$. We define the classifying bisimplicial
smooth manifold as the quotient
$$B_{\bullet \bullet}=U_{\bullet \bullet}/G.$$
We get a principal $G$--bundle $U_{\bullet \bullet}\rightarrow B_{\bullet
\bullet}$, the universal principle $G$-bundle and $\|B_{\bullet \bullet}\|$ is homotopy equivalent to the classifying space $BG$ of $G$.

We define now the connection $\theta_{U_{\bullet
\bullet}}\in\Omega^1_{dR}(U_{\bullet \bullet}; {\mathfrak g})$ on the
universal principal $G$--bundle by the compatible sequence
$\{\theta^{(p)}_{U_{\bullet \bullet}}\}$ defined as
$$\theta^{(p)}_{U_{\bullet \bullet}}=\sum_{j=0}^{p} t_j pr^*_j(\theta)\in \Omega^1_{dR}(\Delta^p\times U_{\bullet p}; {\mathfrak g})$$
where $(t_0, \ldots, t_p)$ are the barycentric coordinates of $\Delta^p$ and
$pr_j\co U_{\bullet p}\rightarrow U_{\bullet 0}$ the canonical projections.

The canonical isomorphism of simplicial manifolds $P_{\bullet}\rightarrow U_{\bullet 0}$ gives a $G$--equivariant map
$$\Psi\co P_{\bullet}\rightarrow U_{\bullet \bullet}$$
and induces a map
$$\psi\co X_{\bullet}=P_{\bullet}/G\rightarrow B_{\bullet \bullet}=
U_{\bullet \bullet}/G.$$
such that $(\Psi, \psi)$ pulls back the principal $G$--bundle $P_\bullet$ over $X_{\bullet}$ and the connection $\theta$ as stated in the theorem.
\end{proof}

\section{Multiplicative cohomology and differential characters}\label{MH}

We will now define general versions of Karoubi's multiplicative cohomology
and Cheeger--Simons differential characters for smooth simplicial manifolds
with respect to any given filtration of the simplicial de Rham complex.
As a special case with respect to the `filtration b\^{e}te' we will
recover the group of Cheeger--Simons differential characters for
smooth simplicial manifolds as introduced by Dupont, Hain and Zucker~\cite{DHZ}.

In general, for a given complex $C^*$ let $\sigma_{\geq p}C^*$
denote the filtration via truncation in degrees below $p$ and
similarly $\sigma_{<p}C^*$ denotes truncation of $C^*$ in degrees
greater or equal $p$. Let us first consider the special case of
Deligne's `filtration b\^ete' \cite{De} for the simplicial de Rham
complex $\Omega_{dR}^*(X_{\bullet})$ of a simplicial manifold
$X_{\bullet}$. The `filtration b\^ete' $\sigma=\{\sigma_{\geq
p}\Omega^*_{dR}(X_{\bullet})\}$ is given as truncation in degrees
below $p$
$$\sigma_{\geq p}\Omega_{dR}^j(X_{\bullet}) =
\begin{cases} 0 & j<p, \\ \Omega^j_{dR}(X_{\bullet}) & j\geq p \end{cases}$$
We define the group of Cheeger--Simons differential characters as follows:

\begin{defn}[See Dupont--Hain--Zucker \cite{DHZ}.]
Let $X_{\bullet}$ be a simplicial smooth manifold and $\Lambda$ be a subgroup of
$\C$.
The group of (mod $\Lambda$) differential characters of degree $k$ of
$X_{\bullet}$ is
given by
$$\hat{H}^{k-1}(X_{\bullet}; \C/\Lambda)= H^k(\cone(\sigma_{\geq k}
\Omega^*_{dR}(X_{\bullet})\rightarrow C^*(X_{\bullet};
\C/\Lambda))).$$
\end{defn}

Now let ${\mathcal F}=\{F^r\Omega^*_{dR}(X_{\bullet})\}$ be any
given filtration of the simplicial de Rham complex. We define the
multiplicative cohomology groups of $X_{\bullet}$ with respect to
${\mathcal F}$ as follows:

\begin{defn}
Let $X_{\bullet}$ be a simplicial smooth manifold, $\Lambda$ be a
subgroup of $\C$ and ${\mathcal
F}=\{F^r\Omega^*_{dR}(X_{\bullet})\}$ be a filtration of
$\Omega^*_{dR}(X_{\bullet})$. The multiplicative cohomology groups
of $X_{\bullet}$ associated to the filtration ${\mathcal F}$ are
given by
$$MH^{2r}_{n}(X_{\bullet}; \Lambda; {\mathcal F})=H^{2r-n}(\cone
(C^*(X_{\bullet}; \Lambda)\oplus
F^r\Omega^*_{dR}(X_{\bullet})\rightarrow C^*(X_{\bullet}; \C))).$$
\end{defn}

In order to be be able to introduce secondary characteristic classes
for connections whose curvature and characteristic forms lie in a
filtration of the simplicial de Rham complex we introduce a more
general version of differential characters associated to any given
filtration. For smooth manifolds these invariants were studied
systematically by the first author in \cite{F}.

\begin{defn}
Let $X_{\bullet}$ be a simplicial smooth manifold, $\Lambda$ be a
subgroup of $\C$ and ${\mathcal
F}=\{F^r\Omega^*_{dR}(X_{\bullet})\}$ be a filtration of
$\Omega^*_{dR}(X_{\bullet})$. The groups of differential characters
(mod $\Lambda$) of degree $k$ of $X_{\bullet}$ associated to the
filtration ${\mathcal F}$ are given by
$$\hat{H}^{k-1}_r(X_{\bullet}; \C/\Lambda; {\mathcal F})=
H^k(\cone(\sigma_{\geq
k}F^r\Omega^*_{dR}(X_{\bullet})\rightarrow C^*(X_{\bullet};
\C/\Lambda))).$$
\end{defn}

The truncation in degrees below $k$ of a complex which is already truncated in
degrees below $k$ leaves it unchanged, hence if ${\mathcal F}$ is Deligne's
`filtration b\^ete' of $\Omega^*(X_{\bullet})$, we recover the ordinary groups
of differential characters of $X_{\bullet}$ as in Definition 2.1.

We have the following main theorem generalizing \cite[Theorem 2.3]{F}.

\begin{thm}
Let $X_{\bullet}$ be a simplicial smooth manifold, $\Lambda$ be a
subgroup of $\C$ and ${\mathcal
F}=\{F^r\Omega^*_{dR}(X_{\bullet})\}$ be a filtration of
$\Omega^*_{dR}(X_{\bullet})$. There exists a surjective map
$$\Xi\co\hat{H}^{2r-n-1}_r(X_{\bullet}; \C/\Lambda; {\mathcal F})\rightarrow
MH^{2r}_n(X_{\bullet};\Lambda; {\mathcal F})$$ whose kernel is the
group of forms in $F^r\Omega_{dR}^{2r-n-1}(X_{\bullet})$ modulo
those forms that are closed and whose complex cohomology class is
the image of a class in $H^*(X_{\bullet};\Lambda)$.
\end{thm}

\begin{proof} Let ${\mathcal A}(F^r)$ and ${\mathcal B}(F^r)$ denote
the cone complexes used in the definition of the groups of differential
characters and multiplicative cohomology associated to the filtration
${\mathcal F}$, that is,
\begin{align*}
{\mathcal A}(F^r)&=
  \cone(\sigma_{\geq k}F^r\Omega^*_{dR}(X_{\bullet})
    \rightarrow C^*(X_{\bullet}; \C/\Lambda)) \\
{\mathcal B}(F^r)&=\cone(C^*(X_{\bullet}; \Lambda)\oplus
F^r\Omega^*_{dR}(X_{\bullet})\rightarrow C^*(X_{\bullet}; \C))
\end{align*}
There is a quasi-isomorphism between the cone complexes
\begin{align*}
&\cone(\sigma_{\geq k}F^r\Omega^*_{dR}(X_{\bullet})\rightarrow
C^*(X_{\bullet}; \C/\Lambda)) \\
&\cone(C^*(X_{\bullet};
  \Lambda)\oplus\sigma_{\geq k}F^r\Omega^*_{dR}(X_{\bullet})\rightarrow
  C^*(X_{\bullet}; \C))
\end{align*}
and we get a short exact sequence of complexes
$$0\rightarrow {\mathcal A}(F^r)\rightarrow {\mathcal B}(F^r)\rightarrow
\sigma_{<k}F^r\Omega^*_{dR}(X_{\bullet})\rightarrow 0$$ where
$\sigma_{<k}$ denotes truncation in degrees greater or equal to $k$.
The statement follows now from the long exact sequence in cohomology
associated to this short exact sequence of complexes, because for
$k=2r-n$ the cohomology group
$$H^{2r-n}(\sigma_{<2r-n}F^r\Omega^*_{dR}(X_{\bullet}))$$ is trivial.
\end{proof}

We can identify the classical Cheeger--Simons differential characters
with multiplicative cohomology groups as follows

\begin{cor}\label{cor2.5}
Let $X_{\bullet}$ be a simplicial smooth manifold and $\Lambda$ be a subgroup of
$\C$. There is an isomorphism
$$\hat{H}^{r-1}(X_{\bullet}; \C/\Lambda)\cong
MH^{2r}_r(X_{\bullet};\Lambda; \sigma)$$
\end{cor}

\begin{proof} 
This is a direct consequence of Theorem 2.4 in the case when $n=r$ and the filtration
${\mathcal F}$ is Deligne's `filtration b\^ete' using the definitions and the quasi-isomorphism
of complexes mentioned in the proof of Theorem 2.4. above 
\end{proof}

In Karoubi's original approach \cite{K2,K3} towards multiplicative
cohomology and differential
characters for a smooth manifold $M$  the complex of modified singular cochains
$$\tilde{C}^*(M;\Z)=\cone(\Omega_{dR}^*(M)\times S^*(M;\Z)\rightarrow S^*(M;\C))$$
is used instead. However this complex is chain homotopy equivalent to the
usual complex of (smooth) singular cochains $S^*(M;\Z)$ of $M$. Again, also in
the more general case of a simplicial smooth manifold $X_{\bullet}$
we can define the complex of modified compatible cochains
${\tilde C}^*(X_{\bullet}; \Z)$ as
$${\tilde C}^*(X_{\bullet};\Z)=\cone(\Omega^*_{dR}(X_{\bullet})\times
C^*(X_{\bullet};\Z)\rightarrow C^*(X_{\bullet};\C)).$$
and proceed as in \cite{K3} or \cite{F} for the definition of multiplicative
cohomology.
But we  can then show that the resulting complex using modified cochains
is quasi-isomorphic to the compatible cochain complex $C^*(X_{\bullet}; \Z)$
used here and the resulting cohomology groups are isomorphic to the ones defined above.

As in the manifold case, it can be shown that the multiplicative cohomology
groups fit in the following long exact sequence (compare \cite{K3})
\begin{multline*}
\cdots\longrightarrow H^{2r-n-1}(\|X_{\bullet}\|;\Lambda)\longrightarrow
  H^{2r-n-1}(\Omega^*_{dR}(X_{\bullet}/F^r\Omega^*_{dR}(X_{\bullet})))\\
\longrightarrow  MH^{2r}_{n}(X_{\bullet}; \Lambda; {\mathcal F})
  \longrightarrow \cdots
\end{multline*}
The groups of differential characters fit also in short exact sequences
analogous to the ones in Cheeger--Simons \cite{CS}, which are again a
special case of the one above.

\begin{rem}
There are several equivalent conventions for the cone of a map of complexes
$f^*\co A^*\rightarrow B^*$. Throughout this paper we will
use the following: $\cone(f^*)^n=A^n\oplus B^{n-1}$ with
differential given by $\d(a,b)=(\d_A a, (-1)^{n+1}f^n(a)+\d_B b)$, where
$a\in A^n$, $b\in B^{n-1}$ and $\d_A,\d_B$ are the differentials in the
complexes $A^*,B^*$ respectively.
\end{rem}

We will discuss some applications to specific examples of simplicial
smooth manifolds. In order to deal with them in a unified way, we
briefly recall the notion of a nerve for a topological category (see
Segal \cite{S} or Dupont \cite{D1,D2}).

Let $\cc$ be a topological category, that is, a small category such
that the set of objects $Ob(\cc)$ and the set of morphisms
$\Mor(\cc)$ are both topological spaces such that
\begin{itemize}
\item[(i)] the source and target maps
$$s,t\co \Mor(\cc)\rightarrow Ob(\cc)$$
are continuous maps.
\item[(ii)] composition of arrows is continuous, that is, if
$$\Mor(\cc)^{\circ}\subseteq \Mor(\cc)\times \Mor(\cc)$$
is the set of pairs $(f, f')$ with $s(f)=t(f')$, the
composition map $\Mor(\cc)^{\circ}\rightarrow \Mor(\cc)$
is a continuous map.
\end{itemize}

Associated to a topological category is a simplicial space $\NN(\cc)_{\bullet}=\{\NN(\cc)_n\}$, the nerve of the category $\cc$. We have\\
$$\NN(\cc)_0=Ob(\cc),\quad
\NN(\cc)_1=\Mor(\cc),\quad
\NN(\cc)_2=\Mor(\cc)^{\circ}$$
and in general
$$\NN(\cc)_n\subseteq \Mor(\cc)\times\cdots\times\Mor(\cc)\quad
  (n \text{ times})$$
is the subset of composable strings of morphisms
$$\bullet\stackrel{f_1}\leftarrow\bullet\stackrel{f_2}\leftarrow\bullet\leftarrow\bullet \cdots \bullet\stackrel{f_n}\leftarrow\bullet,$$
that is, an $n$--tuple $(f_1, f_2, \ldots, f_n)\in \NN(\cc)_n$ if and only if
$s(f_i)=t(f_{i+1})$ for all $1\leq i \leq n-1$.

The face maps $\epsi_i\co \NN(\cc)_n \rightarrow \NN(\cc)_{n-1}$ are given
as
$$\epsi_i(f_1, f_2, \ldots, f_n)= \begin{cases} (f_2, \dots, f_n) & i=0,
\\ (f_1, \ldots, f_i\circ f_{i+1}, \ldots, f_n) & 0<i<n, \\
(f_1, \ldots, f_{n-1}) & i=n.\end{cases}$$
The degeneracy maps $\eta_i\co \NN(\cc)_n \rightarrow \NN(\cc)_{n-1}$ are given as
$$\eta_i(f_1, \ldots, f_n)=(f_1, \ldots, f_{i-1}, \id, f_i, \ldots, f_n),
\quad 0\leq i\leq n.$$

The nerve $\NN$ is a functor from the category of topological categories and continuous functors to the category of simplicial spaces.

\subsubsection*{Classifying spaces of Lie groups}
Let $G$ be a Lie group viewed as a
topological category with one object, that is,
$$Ob(G)=*, \quad  \Mor(G)=G.$$
Furthermore let $\bar{G}$ be the topological category defined as
$$Ob(\bar{G})=G, \quad  \Mor(\bar{G})=G\times G.$$
There is an obvious functor
$$\gamma\co \bar{G}\rightarrow G, \quad  \gamma(g_0, g_1)=g_0g_1^{-1}$$
inducing a map
$$\gamma\co \NN(\bar{G})_{\bullet}\rightarrow \NN(G)_{\bullet}, \quad
\gamma(g_0, \ldots, g_n)=(g_0g_1^{-1}, \ldots, g_{n-1}g_n^{-1})$$
between simplicial smooth manifolds and applying the fat realization functor
gives the universal principal $G$--bundle
$$\gamma_G\co EG\rightarrow BG.$$
Using the simplicial smooth manifold $\NN(G)_{\bullet}$ we can now define

\begin{defn}
Let $G$ be a Lie group, $\Lambda$ a subgroup of $\C$ and let
${\mathcal F}=\{F^r\Omega^*_{dR}(\NN(G)_{\bullet}\}$ be a filtration
of $\Omega^*_{dR}(\NN(G)_{\bullet})$. The multiplicative cohomology
groups of $BG$ associated to the filtration ${\mathcal F}$ are
defined as
$$MH^{2r}_n(BG, \Lambda, {\mathcal F})=MH^{2r}_n(\NN(G)_{\bullet}, \Lambda, {\mathcal F})$$
and the group of differential characters as
$$\hat{H}^{k-1}_r(BG, \C/\Lambda, {\mathcal F})=\hat{H}^{k-1}_r(\NN(G)_{\bullet}, \C/\Lambda, {\mathcal F}).$$
\end{defn}

As in the general case we get the identification from \fullref{cor2.5}
in the case of the `filtration b\^ete' $\sigma$
$$\hat{H}^{r-1}(BG, \C/\Lambda)\cong MH^{2r}_{r}(BG, \Lambda, \sigma).$$
for the classical Cheeger--Simons differential characters.
These invariants were studied in the case $G=GL_n(\C)$ already by Dupont,
Hain and Zucker \cite{DHZ}.

We can generalize this situation much further in the following way.

\subsubsection*{Classifying spaces of Lie groupoids}
Let ${\mathcal G}\co X_1 \to X_0$ be a Lie groupoid, that is, both the set of objects
$X_0$ and the set of morphisms $X_1$ are ${\mathcal
C}^{\infty}$--manifolds and all structure maps are smooth and the
source and target maps are both smooth submersions.

As in the example above we can apply the nerve functor to the category ${\mathcal G}$ and we get again a simplicial smooth manifold $X_{\bullet}=\NN({\mathcal G})_{\bullet}$ where
$$X_n=X_1\times_{X_0}X_1 \times\cdots\times_{X_0}X_1\quad
(n \text{ factors}).$$
Let $B{\mathcal G}$ be the classifying space of ${\mathcal G}$, that is, the
fat realization of the nerve $B{\mathcal G}=\|\NN({\mathcal
G}_{\bullet})\|$. We define

\begin{defn}
Let ${\mathcal G} \co X_1\to X_0$ be a Lie groupoid, $\Lambda$ a
subgroup of $\C$ and let ${\mathcal
F}=\{F^r\Omega^*_{dR}(\NN({\mathcal G})_{\bullet}\}$ be a filtration
of $\Omega^*_{dR}(\NN({\mathcal G})_{\bullet})$. The multiplicative
cohomology groups of $B{\mathcal G}$ associated to the filtration
${\mathcal F}$ are defined as
$$MH^{2r}_n(B{\mathcal G}, \Lambda, {\mathcal F})=MH^{2r}_n(\NN({\mathcal G})_{\bullet}, \Lambda, {\mathcal F})$$
and the group of differential characters as
$$\hat{H}^{k-1}_r(B{\mathcal G}, \C/\Lambda, {\mathcal F})=\hat{H}^{k-1}_r(\NN({\mathcal G})_{\bullet}, \C/\Lambda, {\mathcal F}).$$
\end{defn}

\subsubsection*{Actions of Lie groups on smooth manifolds}
Let $X$ be a ${\mathcal C}^{\infty}$--manifold and $G$ a Lie group which acts smoothly from the left on $X$.
We have a Lie groupoid
$${\mathcal G}\co G\times X\to X$$
with source map $s\co G\times X \rightarrow X, \; s(g, x)=x$, target
map $t\co G\times X\rightarrow X, \; t(g, x)=gx$ and composition map
$$m\co (G\times X)\times_{X} (G\times X) \rightarrow G\times X, (g, hx)(h,x)=(gh,x).$$
This Lie groupoid was studied in detail by Getzler \cite{G} in order
to define an equivariant version of the classical Chern character.
Applying the nerve functor again gives a simplicial manifold, the
homotopy quotient, which allows us to define equivariant versions of
the multiplicative cohomology invariants

\begin{defn}
Let $G$ be a Lie group, acting smoothly on a smooth manifold $X$,
$\Lambda$ a subgroup of $\C$ and let ${\mathcal
F}=\{F^r\Omega^*_{dR}(\NN(G\times X\to X)_{\bullet}\}$ be a
filtration of $\Omega^*_{dR}(\NN(G\times X\to X)_{\bullet})$. The
equivariant multiplicative cohomology groups of $X$ associated to
the filtration ${\mathcal F}$ are defined as
$$MH_G^{2r, n}(X, \Lambda, {\mathcal F})=MH^{2r}_n(\NN(G\times X\to X)_{\bullet}, \Lambda, {\mathcal F})$$
and the group of equivariant differential characters as
$$\hat{H}_G^{k-1, r}(X, \C/\Lambda, {\mathcal F})=\hat{H}^{k-1}_r(\NN(G\times X\to X)_{\bullet}, \C/\Lambda, {\mathcal F}).$$
\end{defn}

We will study secondary theories for classifying spaces of Lie
groupoids and Lie groups in more detail in further papers in view of
applications to foliations, differentiable orbifolds and
differentiable stacks. Equivariant differential characters for
orbifolds of type $[M/G]$ for a smooth manifold $M$ with smooth
action of a Lie group $G$ with finite stabilizers were constructed
and studied systematically by Lupercio and Uribe \cite{LU}. Their
approach follows closely the modified definition of Cheeger--Simons
cohomology due to Hopkins and Singer \cite{HS}. It would be
interesting to study the relation of these invariants with the ones
defined here, especially for different filtrations of the de Rham
complex. Chern--Weil theory for principal $G$--bundles over a Lie
groupoid was systematically analyzed by Laurent-Gengoux, Tu and Xu
\cite{LTX}. This framework can be applied to differentiable stacks
using the general de Rham cohomology of differentiable stacks as
developed by Behrend \cite{B}. The framework developed here allows
the definition of multiplicative cohomology groups and groups of
differential characters for arbitrary filtrations of the de Rham
complex of a differentiable stack, which will be the topic of a
sequel to this paper.

\section{Multiplicative bundles and multiplicative K--theory}\label{MK}

Let $G$ be a Lie group and $\theta_0,\ldots,\theta_q$ be connections on the
principal $G$--bundle
$\pi_{\bullet}\co P_\bullet\rightarrow X_\bullet$,  that is,
$$\theta_j\in \Omega^1_{dR}(P_\bullet;{\mathfrak g})$$
such that for all $p$ and all $0\leq j\leq q$
$$\theta^{(p)}_j\in \Omega^1_{dR}(\Delta^p\times P_p;{\mathfrak g})$$
that is, the restrictions $\theta^{(p)}_j$ are connections on the bundle
$$\Delta^p\times P_p\rightarrow \Delta^p\times X_p.$$
Fix $q$ and let $\Delta^q$ be the standard simplex
in $\R^{q+1}$ parameterized by coordinates $(s_0,\ldots,s_q)$.

\begin{lem}
The form $\sum_{j=0}^q\theta_js_j$
defines a (partial) connection on the pullback
bundle $\pi^*P_\bullet\rightarrow  X_\bullet\times \Delta^q$
where $\pi\co X_\bullet\times \Delta^q \rightarrow X_\bullet$ is
the projection.
\end{lem}

\begin{proof}
For each $m$ the sum $(\sum_{j=0}^q\theta_js_j)^{(m)}=
\sum_{j=0}^q\theta_j^{(m)}s_j$
is a connection on the bundle
$$\Delta^m\times P_m\times\Delta^q\rightarrow  \Delta^m\times X_m\times\Delta^q.$$
We have to verify that the compatibility conditions hold.
The strict simplicial structure on $X_\bullet \times \Delta^q$ is given by
the maps $\epsi_i^\prime=\epsi_i\times \id_{\Delta^q}$ for all $i$,
where  $\epsi_i$  is the map given by the
strict simplicial structure on $X_\bullet$.
We have
$$(\epsi^i\times
\id_{X_m\times\Delta^q})^*(\sum_{j=0}^q\theta^{(m)}s_j)=\sum_{j=0}^q(\epsi^i\times
\id_{X_m})^*\theta_j^{(m)}s_j$$
since the forms $\theta_j^{(m)}s_j$ are in
$$\Omega^1_{dR}(\Delta^m\times P_m; {\mathfrak g})\otimes \Omega_{dR}^0(\Delta^q)\subset \Omega^1_{dR}(\Delta^m\times P_m\times\Delta^q).$$
Now, since the $\theta_j$
satisfy the compatibility conditions we have
$$\sum_{j=0}^q(\epsi^i \times \id_{X_m})^*\theta_j^{(m)}s_j=\sum_{j=0}^q(\id_{\Delta^{m-1}}\times \epsi_i)^*\theta_j^{(m-1)}s_j.$$
As before we have
$$\sum_{j=0}^q(\id_{\Delta^{m-1}}\times \epsi_i)^*\theta_j^{(m-1)}s_j=(\id_{\Delta^{m-1}}\times \epsi_i^\prime)^*(\sum_{j=0}^q\theta^{(m-1)}s_j),$$
 which proves the lemma.
\end{proof}

Given  an invariant polynomial $\Phi$ of degree $k$, we denote by
$$\tilde{\Theta}_q(\Phi;\theta_0,\ldots,\theta_q)\in\Omega^{2k}_{dR}( X_\bullet\times\Delta^q)$$
the characteristic form (on $X_\bullet$) associated to $\Phi$ for the (curvature of the)
connection $\sum_{j=0}^q\theta_js_j$.
When $\Phi$ is understood, we will omit
it from the notation for the above form.
The closed form  $\tilde{\Theta}_q(\theta_0,\ldots,\theta_q)$ is
a family of compatible closed forms
$$\tilde{\Theta}_q^{(m)}(\theta_0,\ldots,\theta_q) \in \Omega_{dR}^{2k}(\Delta^m\times X_m\times\Delta^q).$$
We define a form
$\Theta_q(\theta_0,\ldots,\theta_q)\in\Omega_{dR}^{2k-q}(X_\bullet)$
by
$$\Theta_q(\theta_0,\ldots,\theta_q)=\int_{\Delta^q}\tilde{\Theta}_q(\theta_0,\ldots,\theta_q),$$
that is, $\Theta_q(\theta_0,\ldots,\theta_q)$ is the family of forms
$$\Theta_q^{(m)}(\theta_0,\ldots,\theta_q)=\int_{\Delta^q}\tilde{\Theta}_q^{(m)}(\theta_0,\ldots,\theta_q).$$
These forms satisfy the compatibility conditions since the diagram
\[\xymatrix{
\Omega_{dR}^*(\Delta^m\times
X_m\times\Delta^q)\ar[r]^{\int_{\Delta^q}}\ar[d]_{(\epsi^i\times\id_{X_m\times
\Delta^q})^*}&
\Omega_{dR}^*(\Delta^m\times X_m)\ar[d]^{(\epsi^i\times\id_{X_m})^*}\\
\Omega_{dR}^*(\Delta^{m-1}\times
X_m\times\Delta^q)\ar[r]^{\int_{\Delta^q}}&
\Omega_{dR}^*(\Delta^{m-1}\times X_m)\\
\Omega_{dR}^*(\Delta^{m-1}\times
X_{m-1}\times\Delta^q)\ar[r]^{\int_{\Delta^q}}\ar[u]^{(\id_{\Delta^{m-1}}\times
\epsi_i^\prime)^*}& \Omega_{dR}^*(\Delta^{m-1}\times
X_{m-1})\ar[u]_{(\id_{\Delta^{m-1}}\times \epsi_i)^*}}\] commutes and
the forms $\tilde{\Theta}_q^{(m)}(\theta_0,\ldots,\theta_q) \in
\Omega^{2k}_{dR}(\Delta^m\times X_m\times\Delta^q)$ are compatible.

If we denote by $t$ the variables on the simplices $\Delta^p$, by
$x$ the variables on the manifolds $X_p$ and by $s$ the variables on
the simplex $\Delta^q$, then we can write with an obvious notation
the differential of the complex $\Omega^*_{dR}(\Delta^q\times
X_\bullet)$ as $\d=\d_s + \d_{t,x}$, where $\d_{t,x}$ is the
differential of the complex $\Omega^*_{dR}(X_\bullet)$. Since
$\tilde{\Theta}_q(\theta_0,\ldots,\theta_q)$ is closed, we have
\[\d_{t,x}\tilde{\Theta}_q(\theta_0,\ldots,\theta_q)=-\d_s\tilde{\Theta}_q(\theta_0,\ldots,\theta_q).\]
Then we have
\begin{multline*}
\d_{t,x}\Theta_q(\theta_0,\ldots,\theta_q)
  =\d_{t,x}\int_{\Delta^q}\tilde{\Theta}_q(\theta_0,\ldots,\theta_q) \\
=\int_{\Delta^q}\d_{t,x}\tilde{\Theta}_q(\theta_0,\ldots,\theta_q)
  =-\int_{\Delta^q}\d_s\tilde{\Theta}_q(\theta_0,\ldots,\theta_q).
\end{multline*}
By Stokes theorem the last integral is equal to
$-\int_{\partial\Delta^q}\tilde{\Theta}_q(\theta_0,\ldots,\theta_q)$, so we have
proven the analogue of \cite[Theorem 3.3]{K2}.

\begin{prop}\label{st}In the complex $\Omega^*_{dR}(X_\bullet)$ we have
\[\d\Theta_q(\theta_0,\ldots,\theta_q)=-\sum_{i=0}^q(-1)^i\Theta_{q-1}(\theta_0,\ldots,\hat{\theta_i},\ldots\theta_q).\]
\end{prop}

In particular, for $q=1$ we have that given any two connections on $P_\bullet$, $\theta_0$ and $\theta_1$,
and an invariant polynomial $\Phi$, we can write in a canonical way
\[\Phi(\theta_1)-\Phi(\theta_0)=\d\Theta_1(\Phi;\theta_0,\theta_1).\]
In the sequel it will be convenient to consider formal series of
invariant polynomials (like the  total Chern class for example),
which will then give under the Chern--Weil construction formal sums
of differential forms. We now describe a notation (the same as in Karoubi 
\cite{K2,K3}) to write formulae in this setting in a compact
way. We can write a formal series of invariant polynomials $\Phi$ as
a sum $\sum_r\Phi_r$ with $\Phi_r$ a homogeneous polynomial of
degree $r$. Let ${\mathcal F}=\{F^r\Omega^*_{dR}(X_{\bullet})\}$ be
a filtration of the de Rham complex of $X_\bullet$ and
$$\omega=\sum_r\omega_r,\quad\eta=\sum_r\eta_r$$
be formal sums of forms in $\Omega_{dR}^*(X_\bullet)$ (note that we
do not require that $\omega_r$ is of degree $r$, actually most of
the times this will not be the case). We will write $\omega=\eta
\mbox{ mod } {\mathcal F}$ if and only if for each $r$ we have
\begin{equation}\label{req}
\omega_r-\eta_r\in F^r\Omega_{dR}^*(X_\bullet).
\end{equation}
We will also write $\omega=\eta \mbox{ mod } \tilde{{\mathcal F}}$
when for each $r$ the above equation is satisfied modulo exact forms.

With this notation, all the constructions and proofs in the sequel will
be formally the same both for the case of an invariant polynomial,
where  we will be dealing with forms, and
for a formal series of invariant polynomials, in which case we will
work with formal sums of forms homogeneous degree by homogeneous degree.
Hence we will not distinguish between the two cases in what follows,
writing just $\Phi$ and $\omega$ also for formal sums.

\begin{defn}
Let $\Phi$ be an invariant polynomial (or a formal series) and
${\mathcal F}=\{F^r\Omega^*_{dR}(X_{\bullet})\}$ a filtration of the
de Rham complex of $X_\bullet$. An $({\mathcal
F},\Phi)$--multiplicative bundle (or just a multiplicative bundle
when ${\mathcal F}$ and $\Phi$ are understood) over $X_\bullet$ is a
triple $(P_\bullet,\theta,\omega)$ where $P_\bullet$ is a principal
$G$--bundle over $X_\bullet$, $\theta$ is a connection on $P_\bullet$
and $\omega$ is a (formal series of) form(s) in
$\Omega^*_{dR}(X_\bullet)$ such that
\[\Phi(\theta)=\d\omega \mbox{ mod }{\mathcal F}.\]
\end{defn}

An isomorphism
$f\co(P_\bullet,\theta,\omega)\rightarrow (P_\bullet^\prime,\theta^\prime,\omega^\prime)$
between two multiplicative bundles is an isomorphism $f$ of the underlying bundles
$P_\bullet,P_\bullet^\prime$ such that
\[\omega^\prime-\omega=\Theta_1(\theta,f^*\theta^\prime) \mbox{ mod }
\tilde{{\mathcal F}}.\]
As in Karoubi \cite{K2}, using \fullref{st} to prove transitivity, it follows that
isomorphism is an equivalence relation on multiplicative bundles, so we can
make the following definition.

\begin{defn} We denote by $MK^\Phi(X_\bullet;{\mathcal F})$ the set
of isomorphism classes of multiplicative bundles, the multiplicative K--theory of $X_{\bullet}$
with respect to $({\mathcal F},\Phi)$.
\end{defn}

As usual we will omit $\Phi$ and ${\mathcal F}$ from the notation when there is no risk of ambiguity.

\section{Characteristic classes for secondary theories}

Let $G$ be a Lie group. Given a principal $G$--bundle on a simplicial
smooth manifold $X_{\bullet}$ with a connection $\theta$, and an
invariant polynomial $\Phi$ of homogeneous degree $k$ (for the case
of a  formal series of invariant polynomials one has just to work
degree by degree as in  \fullref{MK}), we will associate
characteristic classes with values in multiplicative cohomology
groups  and in groups of differential characters of $X_{\bullet}$
associated to any filtration of the simplicial de Rham complex
$\Omega^*_{dR}(X_{\bullet})$. This  will generalize the secondary
characteristic classes introduced by Karoubi in the case of smooth
manifolds \cite{K2}.

Let $X_{\bullet}$ be a simplicial smooth manifold, ${\mathcal
F}=\{F^r\Omega^*_{dR}(X_{\bullet})\}$ a filtration of the de Rham
complex, $G$ a Lie group and $\Gamma=(P_{\bullet}, \theta, \eta)$ a
$({\mathcal F},\Phi)$--multiplicative bundle.

The connection $\theta$ on the principal $G$--bundle $\pi_{\bullet}$
is given as $1$--form $$\theta\in \Omega^1_{dR}(P_{\bullet};
\mathfrak{g})$$ as in \fullref{sim}. The characteristic form of
\fullref{CW}
$$\Phi(\theta)\in \Omega^{2k}_{dR}(X_{\bullet})$$
can also as usual be seen as a family of forms
$$\Phi(\theta^{(n)})\in \Omega^{2k}_{dR}(\Delta^n\times X_n; \mathfrak{g})$$
satisfying the compatibility conditions.

Since $\Gamma$ is a multiplicative bundle
we have
\begin{equation}\label{u}\Phi(\theta)=d\eta + \omega,\end{equation}
where the  forms $\eta$ and $\omega$ are also  compatible sequences
$\eta=\{\eta^{(n)}\}$ and $\omega=\{\omega^{(n)}\}$ of differential
forms with $\omega\in F^r \Omega^{2k}_{dR}(X_{\bullet})$ and
$\eta\in \Omega^{2k-1}_{dR}(X_{\bullet})$.

The connection $\theta$ is the pullback of
a connection $\theta_{U_{\bullet \bullet}}$ on $U_{\bullet \bullet}$
by a map $\Psi$ as  in \fullref{uni}.
Let $\Lambda$ be again a subring of the complex numbers $\C$,
and assume that $\Phi$ corresponds under the Chern--Weil map to
a $\Lambda$--valued cohomology class.

For every $n$ the inclusion $\imath_n\co B_{n\bullet}\rightarrow
B_{\bullet\bullet}$ induces isomorphisms in cohomology since
$\|B_{n\bullet}\|$ is homotopy equivalent to the classifying space
of $G$. For every $n$ we also have that
$\imath_n^*\theta_{U_{\bullet\bullet}}= \theta_{U_{n\bullet}}$.
Since the form $\Phi(\theta_{U_{n\bullet}})$ represents the class of
$\Phi$ by \fullref{CW}, and
$\imath_n^*\Phi(\theta_{U{\bullet\bullet}})=\Phi(\theta_{U_{n\bullet}})$,
we have that the form $\Phi(\theta_{U{\bullet\bullet}})\in
\Omega^{*}_{dR}(B_{\bullet\bullet})$ represents the class of $\Phi$.
Then it follows that there exist a compatible cocycle $c\in
C^{2k}(B_{\bullet\bullet};\Lambda)$ and a compatible cochain $v\in
C^{2k-1}(B_{\bullet\bullet};\C)$  such that we have
 \begin{equation}\label{k}
 \delta v=c-\Phi(\theta_{U{\bullet\bullet}}) ,
  \end{equation}
 (where for simplicity we omit from the notation the maps from $\Lambda$--cochains to complex cochains and the
quasi-isomorphism with the de Rham complex).

Since $\Psi^*$ maps compatible cochains (in the bisimplicial sense) to
compatible chains (in the simplicial sense),
the triple $\xi(\Gamma)=(\Psi^*(c), \omega, \Psi^*(v)+\eta)$ defines a cocycle in
the cone complex
$$\cone(C^*(X_{\bullet}; \Lambda)\oplus F^r\Omega^*_{dR}(X_{\bullet})\rightarrow C^*(X_{\bullet}; \C))$$
and since $\omega$ is a form of degree $2k$ also a cocycle in the cone complex
$$\cone(C^*(X_{\bullet}; \Lambda)\oplus\sigma_{\geq 2k}F^r\Omega^*_{dR}(X_{\bullet})\rightarrow  C^*(X_{\bullet}; \C)).$$
The triple $\xi(\Gamma)$ is a cocycle, because we have
$\delta\Psi^*c=\Psi^*\delta c=0$ since $c$ is a cocycle, by \eqref{u} we have
$\d\omega=\d\Phi(\theta)+ \d^2\eta=0$,
and also
\[
\delta\Psi^*v  =  \Psi^*\delta v
                              =  \Psi^*(c-\Phi(\theta_{U_{\bullet\bullet}}))
                              =  \Psi^*c-\Phi(\theta)
                              =  \Psi^*c-(\omega+ \d\eta)
\]
The class of $\xi(\Gamma)$ is independent of the choices of $c$ and $v$:
If $c^\prime$ and $v^\prime$ are other choices satisfying \eqref{k}
then we must have $c-c^\prime=\delta u$ and $\delta u=\delta(v-v^\prime)$.
Then, since $H^{2k-1}(\|B_{\bullet\bullet}\|;\C)$ is trivial, there exists
a compatible cochain $w$ such that $\delta w=u+(v-v^\prime)$. If $\xi^\prime(\Gamma)$
is the cocycle obtained from the different choice, then $\xi(\Gamma)-\xi^\prime(\Gamma)=
(\Psi^*\delta u,0,\Psi^*(v-v^\prime))=\d(\Psi^*u,0,\Psi^*w)$.

Hence for  $2r-m=2k$ we can define the class of the multiplicative bundle
$(P_{\bullet}, \theta, \eta)$ in the
multiplicative cohomology group $MH^{2r}_m(X_{\bullet}, \Lambda, {\mathcal F})$
 to be  the class of $\xi(\Gamma)$.
Similarly the class of $(P_{\bullet}, \theta, \eta)$ in
$\hat{H}^{2k-1}_r(X_{\bullet}; \C/\Lambda; {\mathcal F})$ is the class
of the triple $\xi(\Gamma)$.

\begin{prop}
The classes constructed above are characteristic classes of elements
of $MK^\Phi(X_\bullet;{\mathcal F})$.
\end{prop}

\begin{proof}
The naturality follows from the construction.
We show that for two isomorphic multiplicative bundles
$\Gamma=(P_{\bullet},\theta,\eta)$  and
$\Gamma^\prime=(P_{\bullet}^\prime,\theta^\prime,\eta^\prime)$ the
cocycles $\xi(\Gamma)$ and $\xi(\Gamma^\prime)$ are  cohomologous.
We can assume $P_\bullet=P_\bullet^\prime$, and write
$\Phi(\theta)=\omega +\d\eta$ and $\Phi(\theta^\prime)=\omega^\prime
+\d\eta^\prime$ with $\omega,\omega^\prime \in F^r
\Omega^{2k}_{dR}(X_{\bullet})$. Since the two multiplicative bundles
are isomorphic we have
\[\eta^\prime-\eta=\Theta_1(\Phi;\theta,\theta^\prime)+\sigma+\d\rho\]
with $\sigma\in F^r \Omega^{2k-1}_{dR}(X_{\bullet})$. It follows
that
$$\omega^\prime-\omega=\d(\Theta_1(\Phi;\theta,\theta^\prime)-(\eta^\prime-\eta))
=-\d(\sigma+\d\rho).$$ Let $\Psi^\prime$ be the map pulling back
$\Gamma^\prime$ given by \fullref{uni}, let $c^\prime,v^\prime$
be  the cochains used in the construction for the characteristic
cycle $\xi(\Gamma^\prime)$. Then
$$\xi(\Gamma^\prime)-\xi(\Gamma)=(\Psi^{\prime *}c^\prime-\Psi^*c,
\omega^\prime-\omega,
\Psi^{\prime *}v^\prime-\Psi^*v+\Theta_1(\Phi;\theta,\theta^\prime)+\sigma
+\d\rho)$$
is cohomologous to the triple $\zeta=(\Psi^{\prime *}c^\prime-\Psi^*c,0,
\Psi^{\prime *}v^\prime-\Psi^*v+\Theta_1(\Phi;\theta,\theta^\prime))$
since the two differ by the coboundary of $(0,-\sigma,\rho)$.
We can choose $c^\prime=c$ and $v^\prime=v+\Theta_1(\Phi;\theta_{U{\bullet\bullet}},\theta_{U{\bullet\bullet}}^\prime)$
(where $\theta_{U{\bullet\bullet}}^\prime$ is the connection pulling back
to $\theta^\prime$ under $\Psi^{\prime *}$ given by \fullref{uni})
 satisfying \eqref{k}, hence we have, using also the naturality of
the first transgression form,
$$\zeta=(\Psi^{\prime *}c-\Psi^*c,0,
\Psi^{\prime *}v-\Psi^*v+\Theta_1(\Phi;\Psi^{\prime *}\theta_{U{\bullet\bullet}},\theta^\prime) +\Theta_1(\Phi;\theta,\theta^\prime)).$$
Using \fullref{st}, we have that
\begin{multline*}
\Theta_1(\Phi;\Psi^*\theta_{U{\bullet\bullet}},\Psi^{\prime *}\theta_{U{\bullet\bullet}}) +\d\Theta_2(\Phi;\Psi^*\theta_{U{\bullet\bullet}},\Psi^{\prime *}\theta_{U{\bullet\bullet}},\theta^\prime)=\\
\Theta_1(\Phi;\Psi^{\prime *}\theta_{U{\bullet\bullet}},\theta^\prime)
+\Theta_1(\Phi;\Psi^*\theta_{U{\bullet\bullet}},\theta^\prime).
\end{multline*}
Since $\Psi^\prime$ and $\Psi$ are homotopic, there is a chain homotopy
$H$ between the  induced cochain maps; using $H$ we can write $\zeta$ as
$$(\delta Hc,0, \delta Hv +H\delta v+ \Theta_1(\Phi;\Psi^*\theta_{U{\bullet\bullet}},\Psi^{\prime *}\theta_{U{\bullet\bullet}}) +\d\Theta_2(\Phi;\Psi^*\theta_{U{\bullet\bullet}},\Psi^{\prime *}\theta_{U{\bullet\bullet}},\theta^\prime)).$$
Then $\zeta$ is cohomologous to  $(\delta Hc,0,  +H\delta v+
\Theta_1(\Phi;\Psi^*\theta_{U{\bullet\bullet}},\Psi^{\prime
*}\theta_{U{\bullet\bullet}}))$, and since the transgression forms
$\Theta_1(\cdot;\cdot,\cdot)$ are compatible with chain homotopies
(see Dupont--Hain--Zucker \cite[Appendix A]{DHZ}), the former cocycle
is cohomologous to
$$(\delta Hc,0,H\delta
v+H\Phi(\theta_{U{\bullet\bullet}}))=\d(Hc,0,0)$$ because
$Hc=H(\delta v +\Phi(\theta_{U{\bullet\bullet}}))$ by \eqref{k}.
\end{proof}

\begin{rem}
The above characteristic classes can be slightly generalized in the
following way. Suppose $\Phi$ and $\Phi^\prime$ are formal sums of
invariant polynomials such that every $({\mathcal
F},\Phi)$--multiplicative bundle is also a $({\mathcal
F},\Phi^\prime)$--multiplicative bundle (the main example we have in
mind is the Chern character $ch$ and the total Chern class $c$).
Then by the same procedure we can construct the classes associated
to $\Phi^\prime$ of elements of $MK^\Phi(X_\bullet;{\mathcal F})$
with values in the multiplicative cohomology groups and in the
groups of differential characters associated to ${\mathcal F}$.
\end{rem}
\bibliographystyle{gtart}
\bibliography{link}

\begin{thebibliography}{}
\providecommand\bibmarginpar{\leavevmode\marginpar}
\def\urlstyle#1{{\tt #1}}

\bibitem{B}
\textbf{K\,A Behrend}, \href{http://dx.doi.org/10.1016/j.aim.2005.05.025}
  {\emph{On the de {R}ham cohomology of differential and algebraic stacks}},
  Adv. Math. 198 (2005) 583--622 \xox{MR}{2183389}

\bibitem{CS}
\textbf{J Cheeger}, \textbf{J Simons}, \emph{Differential characters and
  geometric invariants}, from: ``Geometry and topology (College Park, Md.,
  1983/84)'', Lecture Notes in Math. 1167, Springer, Berlin (1985)  50--80
  \xox{MR}{827262}

\bibitem{CM}
\textbf{M Crainic}, \textbf{I Moerdijk},
  \href{http://dx.doi.org/10.1007/s00208-003-0473-2} {\emph{\v{C}ech--{D}e
  {R}ham theory for leaf spaces of foliations}}, Math. Ann. 328 (2004) 59--85
  \xox{MR}{2030370}

\bibitem{De}
\textbf{P Deligne}, \href{http://www.numdam.org/item?id=PMIHES_1971__40__5_0}
  {\emph{Th\'eorie de {H}odge. {II}}}, Inst. Hautes \'Etudes Sci. Publ. Math.
  (1971) 5--57 \xox{MR}{0498551}

\bibitem{D1}
\textbf{J\,L Dupont}, \href{http://dx.doi.org/10.1016/0040-9383(76)90038-0}
  {\emph{Simplicial de {R}ham cohomology and characteristic classes of flat
  bundles}}, Topology 15 (1976) 233--245 \xox{MR}{0413122}

\bibitem{D2}
\textbf{J\,L Dupont}, \emph{Curvature and characteristic classes}, Lecture
  Notes in Mathematics 640, Springer, Berlin (1978) \xox{MR}{0500997}

\bibitem{DHZ}
\textbf{J Dupont}, \textbf{R Hain}, \textbf{S Zucker}, \emph{Regulators and
  characteristic classes of flat bundles}, from: ``The arithmetic and geometry
  of algebraic cycles (Banff, AB, 1998)'', CRM Proc. Lecture Notes 24, Amer.
  Math. Soc., Providence, RI (2000)  47--92 \xox{MR}{1736876}

\bibitem{DJ}
\textbf{J\,L Dupont}, \textbf{H Just}, \emph{Simplicial currents}, Illinois J.
  Math. 41 (1997) 354--377 \xox{MR}{1458178}

\bibitem{E1}
\textbf{H Esnault}, \href{http://dx.doi.org/10.1016/0040-9383(88)90014-6}
  {\emph{Characteristic classes of flat bundles}}, Topology 27 (1988) 323--352
  \xox{MR}{963635}

\bibitem{E2}
\textbf{H Esnault}, \emph{Algebraic differential characters}, from:
  ``Regulators in analysis, geometry and number theory'', Progr. Math. 171,
  Birkh\"auser, Boston (2000)  89--115 \xox{MR}{1724888}

\bibitem{F}
\textbf{M Felisatti}, \href{http://dx.doi.org/10.1023/A:1007812113242}
  {\emph{Differential characters and multiplicative cohomology}}, $K$--Theory
  18 (1999) 267--276 \xox{MR}{1722798}

\bibitem{G}
\textbf{E Getzler}, \href{http://dx.doi.org/10.1006/aima.1994.1081} {\emph{The
  equivariant {C}hern character for non-compact {L}ie groups}}, Adv. Math. 109
  (1994) 88--107 \xox{MR}{1302758}

\bibitem{HS}
\textbf{M\,J Hopkins}, \textbf{I\,M Singer},
  \href{http://projecteuclid.org/getRecord?id=euclid.jdg/1143642908}
  {\emph{Quadratic functions in geometry, topology, and {M}--theory}}, J.
  Differential Geom. 70 (2005) 329--452 \xox{MR}{2192936}

\bibitem{K1}
\textbf{M Karoubi}, \emph{Homologie cyclique et {$K$}--th\'eorie}, Ast\'erisque
  149 (1987) 147 \xox{MR}{913964}

\bibitem{K2}
\textbf{M Karoubi}, \href{http://dx.doi.org/10.1007/BF00534193}
  {\emph{Th\'eorie g\'en\'erale des classes caract\'eristiques secondaires}},
  $K$--Theory 4 (1990) 55--87 \xox{MR}{1076525}

\bibitem{K3}
\textbf{M Karoubi}, \href{http://dx.doi.org/10.1007/BF00961455} {\emph{Classes
  caract\'eristiques de fibr\'es feuillet\'es, holomorphes ou alg\'ebriques}},
  $K$--Theory 8 (1994) 153--211 \xox{MR}{1273841}

\bibitem{LTX}
\textbf{C Laurent-Gengoux}, \textbf{J-L Tu}, \textbf{P Xu}, \emph{Chern--Weil
  map for principal bundles over groupoids}  \xox{arXiv}{math.DG/0401420}

\bibitem{LU}
\textbf{E Lupercio}, \textbf{B Uribe}, \emph{Differential characters on
  orbifolds and string connections {I}: {G}lobal quotients}, from:
  ``Gromov--Witten theory of spin curves and orbifolds'', Contemp. Math. 403,
  Amer. Math. Soc., Providence, RI (2006)  127--142 \xox{MR}{2234887}
  \xox{arXiv}{math.DG/0311008}

\bibitem{S}
\textbf{G Segal}, \href{http://www.numdam.org/item?id=PMIHES_1968__34__105_0}
  {\emph{Classifying spaces and spectral sequences}}, Inst. Hautes \'Etudes
  Sci. Publ. Math. 34 (1968) 105--112 \xox{MR}{0232393}

\end{thebibliography}
\end{document}